\documentclass[11pt]{amsart}
\usepackage{amsfonts,amsmath,amsthm,amscd,amssymb,latexsym,cite,verbatim,texdraw,floatflt,
caption2}
\RequirePackage{hyperref}
\usepackage[a4paper]{geometry}

\title{Selectors of discrete coarse spaces   }
\author{ Igor  Protasov}

\address{I.Protasov: Taras Shevchenko National University of Kyiv, Department of Computer Science and Cybernetics, Academic Glushkov pr. 4d, 03680 Kyiv, Ukraine}
\email{i.v.protasov@gmail.com}

\begin{document}
\begin{abstract} 
Given a coarse space 
$(X, \mathcal{E})$ with the bornology $\mathcal B$ 
of bounded subsets, we extend the coarse structure 
 $\mathcal E$ from $X\times X$ to the natural coarse structure  on 
 $(\mathcal B \backslash \lbrace \emptyset\rbrace)\times (\mathcal B \backslash \lbrace \emptyset\rbrace)$ and say that a 
 macro-uniform
 mapping $f: (\mathcal B \backslash \lbrace \emptyset\rbrace)\rightarrow X$
(resp. $f: [ X]^2 \rightarrow X$) is a selector (resp. 2-selector) of $(X, \mathcal{E})$ if 
$f(A)\in A$ for each $A\in \mathcal B\setminus \lbrace\emptyset\rbrace$ 
(resp. $A \in [X]^2$). We prove that a discrete coarse space  $(X, \mathcal{E})$ admits a selector if and only if 
 $(X, \mathcal{E})$ admits a 2-selector if and only if 
 there exists a linear order $\leq$ on $X$ such that the  family of  intervals $\lbrace [a, b]: a,b\in X, \ a\leq b \}$ is a base  for the  bornology $\mathcal B$.

\end{abstract}
\maketitle



\section{ Introduction}

The notion of selectors comes  from {\it Topology}.
Let $X$ be a topological space, $exp \ X$ be the set of  all non-empty closed subsets of $X$ endowed with some (initially, the Vietoris) topology,  $\mathcal{F}$ be a 
non-empty subset of $exp \ X$. 
A continuous mapping  $f: \mathcal{F} \rightarrow X$
is called  an $\mathcal{F}$-selector of $X$ if $F(A)\in A$ for each $A\in \mathcal{F}$. 
The question on selectors of topological spaces was studied in a plenty of papers, we mention only \cite{b1}, \cite{b4}, \cite{b6}, \cite{b7}.

Formally, coarse spaces, introduced independently and simultaneously in \cite{b8} and \cite{b13},
can be considered as asymptotic counterparts of uniform spaces. But actually this notion is rooted in {\it Geometry} and {\it Geometric Group Theory}, see [13, Chapter 1] and [5, Chapter 4]. At this point, we need some basic definitions.

Given a set $X$, a family $\mathcal{E}$  of subsets of $X\times X$ is called a
{\it  coarse structure} on $X$ if

\begin{itemize}
\item{} each $E \in \mathcal{E}$  contains the diagonal $\bigtriangleup _{X}:=\{(x,x): x\in X\}$ of $X$;
\vspace{3 mm}

\item{}  if  $E$, $E^{\prime} \in \mathcal{E}$  then  $E \circ E^{\prime} \in \mathcal{E}$  and
$ E^{-1} \in \mathcal{E}$,    where  $E \circ E^{\prime} = \{  (x,y): \exists z\;\; ((x,z) \in E,  \ (z, y)\in E^{\prime})\}$,    $ E^{-1} = \{ (y,x):  (x,y) \in E \}$;
\vspace{3 mm}

\item{} if $E \in \mathcal{E}$ and  $\bigtriangleup_{X}\subseteq E^{\prime}\subseteq E$  then  $E^{\prime} \in \mathcal{E}$.
\end{itemize}

Elements $E\in\mathcal E$ of the coarse structure are called {\em entourages} on $X$.

For $x\in X$  and $E\in \mathcal{E}$ the set $E[x]:= \{ y \in X: (x,y)\in\mathcal{E}\}$ is called the {\it ball of radius  $E$  centered at $x$}.
Since $E=\bigcup_{x\in X}( \{x\}\times E[x]) $, the entourage $E$ is uniquely determined by  the family of balls $\{ E[x]: x\in X\}$.
A subfamily ${\mathcal E} ^\prime \subseteq\mathcal E$ is called a {\em base} of the coarse structure $\mathcal E$ if each set $E\in\mathcal E$ is contained in some $E^\prime \in\mathcal E^\prime$.

The pair $(X, \mathcal{E})$  is called a {\it coarse space}  \cite{b13} or  a {\em ballean} \cite{b8}, \cite{b11}.

In this paper, all coarse spaces
 under consideration are supposed to be
 {\it connected}: for any $x, y \in X$, there is $E\in \mathcal{E}$ such $y\in E[x]$.
A subset  $Y\subseteq  X$  is called {\it bounded} if $Y= E[x]$ for some $E\in \mathcal{E}$,
  and $x\in X$.
The family $\mathcal{B}_{X}$ of all bounded subsets of $X$  is a bornology on $X$.
We recall that a family $\mathcal{B}$  of subsets of a set $X$ is a {\it bornology}
if $\mathcal{B}$ contains the family $[X] ^{<\omega} $  of all finite subsets of $X$
 and $\mathcal{B}$  is closed   under finite unions and taking subsets. A bornology $\mathcal B$ on a set $X$ is called {\em unbounded} if $X\notin\mathcal B$.
A subfamily  $\mathcal B^{\prime}$ of $\mathcal B$ is called a base for $\mathcal B$ if, for each $B \in \mathcal B$, there exists $B^{\prime} \in \mathcal B^{\prime}$ such that $B\subseteq B^{\prime}$.

Each subset $Y\subseteq X$ defines a {\it subspace}  $(Y, \mathcal{E}|_{Y})$  of $(X, \mathcal{E})$,
 where $\mathcal{E}|_{Y}= \{ E \cap (Y\times Y): E \in \mathcal{E}\}$.
A  subspace $(Y, \mathcal{E}|_{Y})$  is called  {\it large} if there exists $E\in \mathcal{E}$
 such that $X= E[Y]$, where $E[Y]=\bigcup _{y\in Y} E[y]$.

Let $(X, \mathcal{E})$, $(X^{\prime}, \mathcal{E}^{\prime})$
 be  coarse spaces. 
 A mapping $f: X \to X^{\prime}$ is called
  {\it  macro-uniform }  if for every $E\in \mathcal{E}$ there
  exists $E^{\prime}\in \mathcal{E}^{\prime}$  such that $f(E(x))\subseteq  E^{\prime}(f(x))$
    for each $x\in X$.
If $f$ is a bijection such that $f$  and $f ^{-1 }$ are macro-uniform, then   $f  $  is called an {\it asymorphism}.
If  $(X, \mathcal{E})$ and  $(X^{\prime}, \mathcal{E}^{\prime})$  contain large  asymorphic  subspaces, then they are called {\it coarsely equivalent.}

For a coarse space $(X,\mathcal{E})$, we denote by 
$X^\flat$ the set of all non-empty bounded subsets of $X$, so  
$(X^\flat = \mathcal{B} \backslash \{\emptyset\})$ 
 and  by $\mathcal{E}^\flat$ the coarse structure on $X^\flat$ with the base $\{ E^\flat : E\in  \mathcal{E}\}$, where
 $$(A,B)\in E^\flat \Leftrightarrow A \subseteq E[B], \ \ B\subseteq E[A],$$
and say that $(X^\flat, \mathcal{E}^\flat )$ is the {\it hyperballean} of 
$(X,\mathcal{E})$. For hyperballeans see \cite{b2},\cite{b3}, \cite{b9}, \cite{b10}.

We say that a macro-uniform mapping $f: X^\flat \longrightarrow X$ (resp. $f: [X]^2 \longrightarrow X$)  is a {\it selector} 
(resp. {\it 2-selector}) of $(X,\mathcal{E})$ if $f(A)\in A$ for each $A\in X^\flat$
(resp. $A\in [X]^2$).
We note that a selector is a macro-uniform retraction of $X^\flat$ to $[X]^1$ identified with $X$.

We recall that a coarse space $(X,\mathcal{E})$ is {\it discrete} if, for each $E\in  \mathcal{E}$, there exists a bounded subset $B$ of $(X,\mathcal{E})$ such that $E[x]=\{x\}$  for each $x\in X\setminus B$. Every bornology $\mathcal{B}$ on a set $X$ defines the discrete coarse space $X_\mathcal{B} = (X,\mathcal{E} _\mathcal{B})$, where  $\mathcal{E}_\mathcal{B}$ is a coarse structure with the base 
$\{ E_B: B\in  \mathcal{B}\}$, $E_B [x]=B$ if $x\in B$ and 
$E_B [x]= \{x\}$
if $x\in X\setminus B$. On the other hand,  every discrete coarse space $(X, \mathcal{E})$ coincides with 
$X_\mathcal{B}$, where $\mathcal{B}$ is the bornology of bounded subsets of $(X, \mathcal{E})$.

Our goal is to characterize discrete coarse spaces which admit selectors. 
After exposition of results, we conclude with some comments and open problems.

\section{ Results}

Let $\leq$ be a linear order on a set $X$. We say that 
$(X, \leq)$ is 

\begin{itemize}
\item{} 
{\it right (left) well-ordered} if every subset $Y$ of $X$ has the minimal (maximal) element;

\vspace{3 mm}

\item{} {\it right (left) bounded} if 
 $X$ has the  maximal (minimal) element;

\vspace{3 mm}

\item{} {\it bounded } if $X$ is left and right bounded.

\end{itemize}

\vspace{3 mm}

Every linear order $\leq$ on $X$ defines the bornology $\mathcal{B}_\leq$ on $X$ such that the family $\{[a,b]: a,b\in X, a\leq b\}$, where 
$[a,b]= \{ x\in X: a\leq x \leq b \}$, is a base for $\mathcal{B}_\leq$. Clearly, $X\in \mathcal{B}_\leq$ if and only if $(X, \leq)$ is bounded. 

We say that a 
bornology $\mathcal{B}$ on a set $X$ {\it has an interval base} if  there exists a linear order $\leq$ on $X$ such $\mathcal{B} = \mathcal{B}_\leq$.

\vspace{7 mm}

{\bf Theorem 1. } 
{\it For a bornology $\mathcal{B}$ on a set $X$ and the discrete coarse space 
 $X_{\mathcal{B}}$,
 the following statements are equivalent
 \vspace{5 mm}
 
 $(i)$  $X_{\mathcal{B}}$ admits a selector; 
  
 \vspace{5 mm}
 
 $(ii)$  $X_{\mathcal{B}}$ admits a 2-selector; 
  
\vspace{5 mm}
 
 $(iii)$  $\mathcal{B}$ has an interval base. 
 \vspace{7 mm}

 Proof. }
 If  $X\in \mathcal{B}$ then we have nothing to prove: every mapping $f: \mathcal{B}\backslash \{\emptyset\}\longrightarrow X$ (resp. $f: [X]^2\longrightarrow X $) such that $f(A)\in A$ is a selector (resp. 2-selector) and we take an arbitrary  linear order $\leq$ on $X$
 such that $(X, \leq)$ is bounded. 
 In what follows, $X\notin \mathcal{B}$ so $X_ \mathcal{B}$ is unbounded. 
 The implication $(i)\Rightarrow (ii)$ is evident.
\vspace{5 mm}

$(ii)\Rightarrow (iii)$ We take a 2-selector $f$ of $X_ \mathcal{B}$ and define a binary relation $\prec$
on $X $  as  follows: $a\prec b$
if and only if either $a=b$ or $f(\{a,b\}) = a$.
\vspace{5 mm}

We use the following key observation 

\vspace{5 mm}
$(\ast)$ {\it for every $B\in \mathcal{B}$, there exists $C\in \mathcal{B}$  such that if $z\in X\setminus C$ then either  $b\prec z$ for each $b\in B$ or $z\prec b$ for each  $b\in B$.}
\vspace{5 mm}

Indeed, we take $C\in \mathcal{B}$ such that 
$B\subseteq C$  and  if  $A, A^\prime \in [X]^2$ and 
$(A, A^\prime) \in E_B ^\flat$ then 
$(f(A), f(A^\prime)) \in E_C$.

\vspace{10 mm}

We take and fix distinct  $l, r \in X$ such that 
$l\prec r$ and use the Zorn's lemma to choose a maximal by inclusion subset $A$ of $X$ such that $A=L\cup R$,
$L\cap R = \emptyset$, $R$ is right well-ordered by $\prec$ with the minimal element $r$, $L$ is left 
well-ordered by $\prec$ with the maximal element $l$
and $x\prec y$ for all $x\in L$, $y\in R$.

By the maximality of $A$ and $(\ast)$, $A$ is unbounded in $X_\mathcal{B}$. 
For  $a,b\in A$,  $a\prec b$, we denote $[a,b]_A = \{ x\in A : a\prec x\prec b\}$.
Applying $(\ast)$ with $B= [a,b]_A$, we see that 
$[a,b]_A$ is bounded in $X_\mathcal{B}$. 

We consider three cases.

\vspace{5 mm}

{\it Case 1:} $L$ and $R$ are  unbounded in 
$X_\mathcal{B}$. We define  some auxiliary mapping $h: X \longrightarrow A$.  For $x\in A$, we put $h(x)=x$.
For $x\in X\backslash A$, we use $(\ast)$  with $B=\{r, x\}$ to find the minimal element $c\in R$ such that $x\prec y$ for each $y\in A$, $c\prec y$. 
If $c\neq r$ then we put $h(x)=c$. 
Otherwise, we use $(\ast)$ to choose the maximal element $d\in L$  such  that 
$y\prec x$ for each 
$y\in L\cup \{r\} $, $y\prec d$. 
We put $h(x)=d$. 

We take arbitrary $a,b\in A$ such that $a\prec l\prec r\prec b$.
If $h(x)\in [a,b]$ then, by the construction of $h$, we have $a\prec x\prec  b$.
Applying $(\ast)$ with $B= [a,b]_A$, we conclude that 
$h^{-1} ([a,b]_A)$ is bounded in $X_\mathcal{B}$.
In particular, $h^{-1}(c)$ is bounded in 
$X_\mathcal{B}$ for each $c\in A$.

Now we are ready to define the desired linear order $\leq$ on $X$.  
If $h(x) \prec h(y)$  and $h(x) \neq h(y)$ then we put $x<y$.
If $c\in R$ then we endow $h^{-1}(c)$ with a right well-order $\leq$. 
If $c\in L$ then we endow  $h^{-1}(c)$ with a left well-order $\leq$.

It remains to verify that the family 
$\{[a,b] : a,b \in X$, $\ a\leq b \}$ is a base for 
$\mathcal{B} $. Let $a,b\in A$ and $a\leq b$. We  have shown that $h^{-1} ([a, b ]_A)\in \mathcal{B}$, hence 
$[a, b ]\in \mathcal{B}$.
If $a,b\in X$  and $a\leq b$ then we take 
$a^\prime \in A$, $b^\prime \in A $
such that $a^\prime < a$, $b< b^\prime$. 
Since $[a^\prime , b^\prime ] \in \mathcal{B}$, we
have $[a, b] \in \mathcal{B}$.
On the other hand, if $Y$ is a bounded subset of  
$X_\mathcal{B}$ then we apply $(\ast)$
 with $B= Y \cup \{l, r\}$
 to find $a\in L$, $b\in R$  such that 
 $h(B) \subseteq [a,b]_A$, hence $B \subseteq [a,b]$.
 
 \vspace{7 mm}

{\it Case 2:} $L$ is bounded and $R$ is unbounded in $X_\mathcal{B}$.
Since $L\in \mathcal{B}$, by $(\ast)$,
the set $C=\{ x\in X: x< y $ for each $y\in R\}$ is bounded in  $X_\mathcal{B}$.
We use arguments  from Case 1 to define $\leq$ on $X\setminus C$.
Then we extend  $\leq$ to $X$ so that $(C, \leq)$ is bounded and $x\prec y$ for all $x\in C$, $y\in X\setminus C$.

  \vspace{7 mm}

{\it Case 3:} $L$ is unbounded and $R$ is bounded in $X_\mathcal{B}$.
Since 
$R\in \mathcal{B}$, by $(\ast)$,
the set 
$D=\{ x\in X: y\prec x $ for 
each $y\in L\}$ is bounded in  $X_\mathcal{B}$.
We use arguments  from Case 1 
to define $\leq$ on $X\setminus D$.
Then we extend 
 $\leq$ to $X$
  so that $(D, \leq)$ is bounded and $x\prec y$ for all $x\in X\setminus D$, $y\in D$.
 
 \vspace{7 mm}
 
 $(iii) \Rightarrow (i)$. We  take a linear order $\prec $ on $ X$ witnessing that $\mathcal{B}$ has an interval base.
 We define a 2-selector $f: [X]^2 \longrightarrow X$  by $f(\{x,y\}) = x $ if and only if $x\prec y$.
 Then we take the linear order $\leq$ on $X$
 defined in the proof $(ii)\Rightarrow (iii)$.
 To define a selector $s$ of $X_\mathcal{B}$,
 we denote $X_l = \{ x\in X: x\leq l\}$, $X_r = \{ x\in X: r\leq x\}$.
 By the construction of $\leq$, $X_l$ is right  well-ordered and $X_r$ is left well-ordered. We take an arbitrary   $Y\in \mathcal{B}\setminus \{\emptyset \}$.
 If $Y\cap X_l\neq \emptyset$ then we take the maximal element $a\in Y\cap X_l$ and put $s(Y)=a$.
 Otherwise, we choose the minimal element $b\in Y\cap X_r$ and put $s(Y)=b$.
 
 To see that $s$ is macro-uniform, we take an interval 
 $[a, b]$ in $(X, \leq)$ and $Y,Z\in \mathcal{B}\setminus \{\emptyset \}$
 such that $Y\setminus [a,b]= Z \setminus [a,b]$,
 $Y\cap [a,b]\neq\emptyset,$  
 $Z\cap [a,b]\neq\emptyset$.
 If  $s(Y)\notin [a,b]$ then  $s(Y) = s(Z)$.
 If  $s(Y)\in [a,b]$  then 
 $s(Z)\in [a,b]$.
 $ \ \  \  \Box $
 \vspace{7 mm}

An ordinal $\alpha$ endowed with the reverse ordering is called the {\it antiordinal} of $\alpha $.

 \vspace{7 mm}
 {\bf Corollary 2. } 
{\it If $X_\mathcal{B}$ has a selector then 
$\mathcal{B}$ has an interval base with respect to some linear order $\leq$ on $X$ such that $(X, \leq)$ is the ordinal sum of an antiordinal and an ordinal.
\vspace{3 mm}

 Proof. } We take the linear order from the proof of Theorem 1 and note that $X_l$ is an antiordinal, $X_r$ is ordinal and $(X, \leq)$ is the ordinal sum of  $X_l$ and $X_r$.
$ \ \  \  \Box $
 \vspace{7 mm}
 
 {\bf Corollary 3. } 
{\it If a bornology $\mathcal{B}$ on a set $X$ has a base linearly ordered by inclusion then the discrete coarse space $X_\mathcal{B}$ admits a selector. 
 
 \vspace{3 mm}

 Proof. } 
 Since $\mathcal{B}$ has a linearly ordered base, we can choose a base $\{ B_\alpha : \alpha < \kappa \}$
 well-ordered by inclusion. 
 We show that $\mathcal{B}$  has an interval base and apply Theorem 1.
 
 For each $\alpha < \kappa$, let $\mathcal{D}_\alpha = B_{\alpha +1} \setminus B_{\alpha}$. We endow  each 
 $\mathcal{D}_\alpha $ 
 with an arbitrary right well-order $\leq$.
 If $x\in \mathcal{D}_\alpha $, 
 $y\in \mathcal{D}_\beta $ and $\alpha< \beta$, we put $x<y$. 
 Then $\mathcal{B}= \mathcal{B}_\leq $. 
 $ \ \  \  \Box $
 \vspace{7 mm}
 
 {\bf Remark 4. } Let $(X, \leq)$ be the ordinal sum of the antiordinal of $\omega$ and the ordinal $\omega_1$.
 Then the interval bornology 
 $\mathcal{B}_\leq $ does not have a linearly ordered base. Indeed, let $X=L\cup R$, $L=\{ l_n : n<\omega\}$,
 $l_n< l_m$ iff $m< n$, $R=\{ r_\alpha : \alpha< \omega_1 \}$, $r_\alpha< r_\beta$ iff $\alpha<\beta$,
 and $l_n < r_\alpha$ for all $n, \alpha$.
 Assuming that $\mathcal{B} _{\leq}$ has a linearly 
 ordered base,  we choose a base $\mathcal{B}^\prime $
 of  $\mathcal{B} _{\leq}$  well-ordered by inclusion and denote $\mathcal{B}^\prime _n =\{A\in \mathcal{B}^\prime : min \ A=l_n \}$.
 By the choice of $R$, there exists $m\in \omega$ such that 
 $\mathcal{B}^\prime _m$ is cofinal in $\mathcal{B} _{\leq}$,  but $l_{m+1}\notin A$ for each 
 $A\in \mathcal{B}^\prime _m$ and we get a contradiction. 
 
 \vspace{7 mm}
 
 {\bf Theorem 5. }{\it Let $(X, \mathcal{E})$ be a coarse space with the bornology  $\mathcal{B}$ of bounded subsets. 
 If  $f$ is a 2-selector of $(X, \mathcal{E})$ then $f$
 is a 
 2-selector of $X_\mathcal{B}$.
 
 \vspace{3 mm}

 Proof. } Let $B\in \mathcal{B}$, $A, A^\prime \in [X]^2$ and $(A, A^\prime )\in E_B ^\flat$. 
 Since $f$ is a 2-selector of $X_\mathcal{E})$, there exists $F\in \mathcal{E}$, $F=F^{-1}$
 such that 
 $(f(A), f(A^\prime ))\in F$.
 
 If $A\cap B=\emptyset$ then  $A=A^\prime$.
 If $A\subseteq B$ then  $A^\prime\in B$, so $(f(A), f(A^\prime ))\in E_B$.
 
 Let  $A=\{b,a\}$, $A^\prime=\{b^\prime,a\}$, $b\in B$,
 $b^\prime\in B$ and $a\in X\setminus B$.
 If $a\in F[\{b,b^\prime\}]$ then 
  $f(A), f(A^\prime )\in F[\{b,b^\prime\}]$.
  If $a\notin F[\{b,b^\prime\}]$
 then either  
$f(A)=f(A^\prime )=a$  or $f(A), f(A^\prime )\in  \{b, b^\prime \}$.

In all considered cases, we have  $(f(A), f(A^\prime ))\in E_{F[B]}$.
 Hence, $f$ is a 2-selector of $X_\mathcal{B}$.
 $ \ \  \  \Box $
 \vspace{7 mm}

 {\bf Remark 6. } Every metric space $(X, d)$ has the natural coarse structure $\mathcal{E}_d$ with the base 
 $\{E_r: r>0 \}$, $ \ E_r=  \{(x,y); d(x,y)\leq r \}$.
 Let $\mathcal{B}$ denotes the bornology of bounded subsets of $(X, \mathcal{E}_d)$.
 By Corollary 3, the discrete coarse space $X_\mathcal{B}$ admits a 2-selector. 
 We show that 
 $(X, \mathcal{E}_d)$ could not admit a 
 2-selector, so the conversion of Theorem 5 does not hold.
 
 Let $X=\mathbb{Z}^2$, 
 $x=(x_1 , x_2)$, $y=(y_1 , y_2)$, 
 $d(x,y)= max \ \{ |x_1 - y_1|, |x_2 - y_2|\}$. 
 We suppose that there exists a 2-selector $f$  of 
 $(X, \mathcal{E}_d)$ and choose a natural number $n$
 such that if $A, A^\prime\in [X]^2$
 and 
 $(A, A^\prime )\in E_1 ^\flat$ then 
 $(f(A), f(A^\prime ))\in E_n$, so
 $d(f(A), f(A^\prime ))\leq n$.
 We denote $S_n=\{ x\in X : d(x,0)=n\}$.
 For $x\in S_n$,  let $A_x = \{x, -x\}$.
 Then we can choose $x,y\in S_n$ such that 
 $d(x,y)=1$, $f(A_x)=x$, $f(A_y)=-y $, but $d(x, -y)>n$.

\section{Comments }



\vspace{7 mm}

1. Let $(X, \mathcal{U})$ be a uniform  space and let $\mathcal{F}_X$ denotes the set of all non-empty closed subsets of $X$
endowed with the Hausdorff-Bourbaki uniformity. 
Given a subset $\mathcal{F}$ of $\mathcal{F}_X$, 
a uniformly continuous mapping 
$f: \mathcal{F} \longrightarrow X$ is called 
an $\mathcal{F}$-selector if $f(A)\in A$ for each $a\in A$. 
If  $\mathcal{F}=[X]^2$
then $f$ is called a 2-selector.

In contrast to the topological case, the problem of uniform selections is much less studied. Almost all known results are concentrated around uniformizations of the Michael's theorem, for references see \cite{b12}.

Given a discrete uniform space, how can one detects whether $X$ admits a 2-selector?
This questions, seems, very difficult even in the case of a countable discrete metric space $X$.
To demonstrate the obstacles for a simple characterization, we consider the following example. 

We take a family $\{ C_n : n<\omega \}$ of pairwise  non-intersecting circles of radius 1 on the Euclidean plane $\mathbb{R}^2$. Then we  inscribe a regular $n$-gon $M_n$ in $C_n$ and denote by $X$ the set of all vertices   of $\{ M_n : n<\omega \}$.
It is  easy  to verify  that $X$ does not admit a 2-selector.

\vspace{7 mm}


2. Given a group $G$ with the identity $e$, we denote by 
$\mathcal{E}_G$ a coarse structure on $G$ with the  base $$\{\{ (x,y) \in G\times G: y\in Fx \}: F\in [G]^{<\omega} , \ \ e\in F \}$$ and say that $(G, \mathcal{E}_G)$ is the {\it finitary  coarse space} of $G$. It should be 
noticed that finitary  coarse spaces of groups 
(in the form of Cayley graphs)
are used in {\it Geometric Group Theory}, see  \cite{b5}.
We note that the bornology of bounded subset of $(G, \mathcal{E}_G)$ is the set $[G]^{<\omega}$. Applying 
Theorem 1 and Theorem 5, we conclude that if  
$(G, \mathcal{E}_G)$
admits a 2-selector
 then 
$G$ must be countable.

\vspace{5 mm}

{\bf Problem 1.} {\it Characterize countable groups  $G$ such that the finitary coarse space $(G,\mathcal{E}_G)$ admits a 2-selector.}

\vspace{7 mm}

3. Every connected graph $\Gamma [\mathcal{V}]$ with the set of vertices $\mathcal{V}$ can be considered as the 
metric space $(\mathcal{V},d)$, where $d$ is the path metric on $\mathcal{V}$.

\vspace{5 mm}

{\bf Problem 2.} {\it Characterize  graphs  $\Gamma [\mathcal{V}]$ such that 
the  coarse space of $(\mathcal{V}, d)$, where $d$ is the path metric on the set of  vertices  $\mathcal{V}$,
admits a 2-selector.}

\vspace{5 mm}

{\bf Acknowledgments.} I thank the referee for critical remarks to the initial version of the paper.


\end{document}